\documentclass[a4paper]{amsart}

\usepackage{
 amssymb, 
 fullpage,
 mathtools, 
 varioref, 
 tikz}
\usetikzlibrary{positioning}
\usetikzlibrary{arrows}
\usepackage[pagebackref]{hyperref}

\theoremstyle{plain}
    \newtheorem{theorem}{Theorem}[section]
    \newtheorem{proposition}[theorem]{Proposition}
    \newtheorem{definition}[theorem]{Definition}
    \newtheorem{conjecture}[theorem]{Conjecture}
\theoremstyle{remark}
    \newtheorem{remark}[theorem]{Remark}

\DeclareMathOperator{\Hom}{Hom}
\DeclareMathOperator{\Fil}{Fil}
\DeclareMathOperator{\GL}{GL}
\DeclareMathOperator{\loc}{loc}
\DeclareMathOperator{\Gal}{Gal}
\DeclareMathOperator{\rank}{rank}

\newcommand{\cF}{\mathcal{F}}
\newcommand{\cK}{\mathcal{K}}
\newcommand{\cN}{\mathcal{N}}
\newcommand{\cO}{\mathcal{O}}
\newcommand{\cP}{\mathcal{P}}

\newcommand{\fm}{\mathfrak{m}}
\newcommand{\fp}{\mathfrak{p}}
\newcommand{\fq}{\mathfrak{q}}

\newcommand{\CC}{\mathbf{C}}
\newcommand{\QQ}{\mathbf{Q}}
\newcommand{\ZZ}{\mathbf{Z}}
\newcommand{\DD}{\mathbf{D}}
\newcommand{\FF}{\mathbf{F}}

\newcommand{\bfc}{\mathbf{c}}
 
\newcommand{\Qp}{{\QQ_p}}
\newcommand{\Zp}{{\ZZ_p}}

\newcommand{\f}{\mathrm{f}}
\newcommand{\dR}{\mathrm{dR}}
\newcommand{\Iw}{\mathrm{Iw}}
\newcommand{\ord}{\mathrm{ord}}

\numberwithin{equation}{section}

\title[Euler systems with local conditions]{Euler systems with local conditions}

\author{David Loeffler}
\address[Loeffler]{Mathematics Institute\\
 Zeeman Building, University of Warwick\\
 Coventry CV4 7AL, UK.\\ 
ORCID identifier: \href{https://orcid.org/0000-0001-9069-1877}{0000-0001-9069-1877}}
\email{d.a.loeffler@warwick.ac.uk}

\author{Sarah Livia Zerbes}
\address[Zerbes]{Department of Mathematics \\
 University College London\\
 Gower Street, London WC1E 6BT, UK.\\
ORCID identifier: \href{https://orcid.org/0000-0001-8650-9622}{0000-0001-8650-9622}}
\email{s.zerbes@ucl.ac.uk}

\thanks{The authors are grateful to acknowledge financial support from the European Research Council (ERC Consolidator Grant ``Euler Systems and the Birch--Swinnerton-Dyer conjecture''), the Royal Society (University Research Fellowship ``L-functions and Iwasawa theory''), and the Leverhulme Trust (Philip Leverhulme Prize PLP-2014-354).}

\begin{document}

\begin{abstract}
 Euler systems are certain compatible families of cohomology classes, which play a key role in studying the arithmetic of Galois representations. We briefly survey the known Euler systems, and recall a standard conjecture of Perrin-Riou predicting what kind of Euler system one should expect for a general Galois representation. Surprisingly, several recent constructions of Euler systems do not seem to fit the predictions of this conjecture, and we formulate a more general conjecture which explains these extra objects. The novel aspect of our conjecture is that it predicts that there should often be Euler systems of several different ranks associated to a given Galois representation, and we describe how we expect these objects to be related.
\end{abstract}

\subjclass[2010]{11R23, 11F80}

\maketitle

\subsubsection*{Acknowledgements} This survey is based on the notes of a lecture given by the second author at the 60th birthday conference for Michael Harris at MSRI, Berkeley, in 2014. It was expanded to a survey article during the ``Euler Systems and Special Cycles'' semester at the Centre Bernoulli, Lausanne, in autumn 2017. We are grateful to MSRI and to the Centre Bernoulli for their hospitality. We would also like to thank the anonymous referee, and the editorial board of the Iwasawa centenary proceedings, for their valuable comments on earlier drafts of this paper.

\section{Cohomology of Galois representations}
  
 The representations of Galois groups of number fields play a central role in number theory. For instance, if $K$ is a number field and $E / K$ is an elliptic curve, one can consider its \emph{Tate module} 
 \[ T_p(E) \coloneqq \varprojlim_n E[p^n], \]
 for a prime $p$; this is a free rank 2 $\Zp$-module, with a continuous action of the group $G_{K} \coloneqq \Gal(\overline{K} / K)$. The representation $T_p(E)$ contains much useful arithmetic data about $E$; for instance, $E$ has good reduction at a prime $\ell \ne p$ if and only if the inertia group $I_\ell$ acts trivially on $T_p(E)$ (the ``Ner\'on--Ogg--Shafarevich criterion''). 
 
 Deeper properties of $E$ are encoded in the (continuous) Galois cohomology groups $H^i(G_K, T_p(E))$, which we shall abbreviate as $H^i(K, T_p(E))$ henceforth. There is a natural injective map, the \emph{Kummer map},
 \[ \kappa: E(K) \otimes \Zp \to H^1(K, T_p(E)), \]
 and many of the deepest results we have concerning the Mordell--Weil groups of elliptic curves -- notably Kolyvagin's theorem that if $E$ is an elliptic curve over $\QQ$ and $\ord_{s = 1} L(E, s) \le 1$, then the Tate--Shafarevich group of $E$ is finite and the Birch--Swinnerton-Dyer conjecture holds for $E$ -- have been proved by studying the image of $E(K)$ in $H^1(K, T_p(E))$, using sophisticated techniques in Galois cohomology. So describing and controlling the cohomology of Galois representations is a deep and fundamental problem.
 
 One of the few tools available for controlling global cohomology groups is the theory of \emph{Euler systems}, and in this article we shall introduce the theory of Euler systems, and formulate a new conjecture predicting what sort of Euler systems one might expect for general Galois representations. 
 
\section{Euler systems}
 \label{sect:ES}
 The definition of an Euler system comes in several slightly different flavours. We shall follow the standard reference, which is \cite{rubin00}. As above, we fix a number field $K$, a prime $p > 2$, and a finite extension $L / \Qp$ with ring of integers $\cO$.

 \begin{definition}
  For an integral ideal $\fm$ of $K$, we write $K(\fm)$ for the maximal extension of $K$ of $p$-power degree contained the ray class field of $K$ modulo $\fm$.
 \end{definition}

 Let $T$ be a finite-rank free $\cO$-module with a continuous action of $G_K$, unramified at almost all primes. We write $T^*(1)$ for the Tate dual $\Hom(T, \cO(1))$, and if $\fq$ is a finite prime of $K$ at which $T$ is unramified, we define a local Euler factor $P_{\fq} \in \cO[X]$ by
 \[
  P_{\fq}(X) \coloneqq \sideset{}{_{\cO}}\det\left(1 - X \operatorname{Frob}_\fq^{-1} : T^*(1)\right).
 \]
 
 We fix an ideal $\cN$ of $K$, divisible by $p$ and by all primes at which $T$ is ramified; and an infinite abelian extension $\cK$ of $K$ which contains $K(\fq)$, for every prime $\fq \nmid \cN$, and the cyclotomic $\ZZ_p$-extension $K_\infty \subset K(p^\infty)$. 
 
 \begin{definition}[{\cite[Definition 2.1.1]{rubin00}}]
  \label{def:ES}
  An \emph{Euler system} for $(T, \cK, \cN)$ is a collection of cohomology classes
  \[ \bfc = \{ \bfc_F \in H^1(F, T): K \subseteq_{\f} F \subset \cK \}, \]
  (where the notation $K \subseteq_{\f} F \subset \cK$ signifies that $F$ runs over the finite extensions of $K$ contained in $\cK$), satisfying the following relation: if $K \subseteq_{\f} F \subseteq_{\f} F' \subset \cK$, then
  \[ \operatorname{cores}_F^{F'}\left(\bfc_{F'}\right) = \left( \prod_{\fq \in \Sigma(F' / F)}  P_\fq\left(\sigma_{\fq}^{-1}\right)\right) \bfc_F \tag{$\star$}\]
  where $\Sigma(F' / F)$ is the set of (finite) primes of $K$ not dividing $\cN$ which ramify in $F'$ but not in $F$, and $\sigma_{\fq}$ is the image of $\operatorname{Frob}_\fq$ in $\Gal(F / K)$.
 \end{definition}
 
 Note that only the local Euler factors at unramified primes appear in the definition; the Euler factors at the bad primes play no direct role. 
 
 \begin{remark}
  As noted in \cite[\S 9.4]{rubin00}, Kolyvagin's Euler system of Heegner points does not actually fit into the definition \ref{def:ES}, since Heegner points are always defined over abelian extensions of an imaginary quadratic field $K$ which are \emph{anticyclotomic} -- one cannot find Heegner points defined over all the fields $K(\fq)$ in such a way that the Euler system norm relations are satisfied. There are other examples of anticyclotomic Euler systems, but we shall not discuss them further in this survey, for reasons of space.
 \end{remark}
 
 The basic function of Euler systems is to bound \emph{Selmer groups}, which are subgroups of $H^1(K, T)$ defined by local conditions.
 
 \begin{definition}\mbox{~}
  \begin{enumerate}
   \item[(i)] If $v \nmid p$ is a (finite) prime of $K$, we define
   \[ H^1_{\f}(K_v, T) \coloneqq \ker\left(H^1(K_v, T) \to H^1(K_v^{\mathrm{nr}}, T \otimes \Qp)\right)\]
   where $K_v^{\mathrm{nr}}$ is the maximal unramified extension of $K_v$.
   \item[(ii)] We define the \emph{relaxed Selmer group}, $H^1_{\mathrm{rel}}(K, T)$, as
   \[
    \left\{ 
     x \in H^1(K, T): 
     \loc_v(x) \in H^1_\f(K_v, T) \text{ for all $v \nmid p$}
    \right\}.
   \]
   \item[(iii)] We define the \emph{strict Selmer group}, $H^1_{\mathrm{str}}(K, T)$, as
   \[
    \left\{ 
     x \in H^1(K, T): 
     \begin{array}{ll}
      \loc_v(x) \in H^1_\f(K_v, T) &\text{for all $v \nmid p$},\\
      \loc_v(x) = 0 &\text{for all $v \mid p$.}
     \end{array}
    \right\}
   \]
  \end{enumerate}
 \end{definition}
 
 \begin{theorem}[{Rubin, cf.~\cite[Theorem 2.2.3]{rubin00}}]
  \label{thm:rubin}
  Suppose $\bfc$ is an Euler system for $(T, \cK, \cN)$, and $\bfc_K$ is non-torsion in $H^1(K, T)$; and suppose that $T$ satisfies a mild ``large image'' hypothesis. Then the group $H^1_{\mathrm{str}}(K, T^*(1))$ is finite.
 \end{theorem}
 
 \begin{remark}
  Rubin states his theorem in a somewhat different form, involving the finiteness of the strict Selmer group of the $p$-torsion representation $T^*(1) \otimes \Qp/\Zp$, but this is equivalent to the above statement.
 \end{remark}
 
 The Poitou--Tate global duality theorem for Galois cohomology, combined with Tate's Euler characteristic formula, shows that the finiteness of the strict Selmer group of $T^*(1)$ implies a bound for the cohomology of $T$. This bound involves the following important numerical invariant:
 
 \begin{definition}
  We define 
  \[ d_-(T) \coloneqq \sum_{\substack{v \mid \infty \\ \text{$v$ real}}} \rank_{\cO}\left( T^{\sigma_v = -1}\right) + \sum_{\substack{v \mid \infty \\ \text{$v$ complex}}} \rank_{\cO}\left(T\right),\]
  where $\sigma_v$ denotes complex conjugation at $v$.
 \end{definition}
 
 Let us write $h^i(-)$ for the rank of the cohomology group $H^i(-)$ as an $\cO$-module.
 
 \begin{proposition}
  \label{prop:h1bound}
  Suppose $H^0(K, T) = H^0(K, T^*(1)) = 0$, and $H^0(K_v, T^*(1)) = 0$ for all primes $v \mid p$. Then we have
  \[ h^1_{\mathrm{rel}}(K, T) \ge d_-(T),\]
  with equality if and only if $H^1_{\mathrm{str}}(K, T^*(1)) = 0$.
 \end{proposition}

 So, by Rubin's theorem, the existence of a non-trivial Euler system forces $H^1_{\mathrm{rel}}(K, T)$ to have the minimal possible rank.
 
 \begin{proof}
  Let $V \coloneqq T \otimes_\cO L$, let $S$ be the set of primes dividing $\cN \infty$, and let $G_{K, S}$ be the Galois group of the maximal extension of $K$ unramified outside $S$. Then, for any $T$ unramified outside $S$, Poitou--Tate duality gives an exact sequence of finite-dimensional $L$-vector spaces
  \begin{gather*}
   0 \to H^1_{\mathrm{rel}}(K, V) \to H^1(G_{K, S}, V) \to \bigoplus_{\substack{v \in S \\ v \nmid p}} H^1_{\mathrm{s}}(K_v, V) 
   \to H^1_{\mathrm{str}}(K, V^*(1))^* \\ \to H^2(G_{K, S}, V) \to \bigoplus_{v \in S} H^2(K_v, V) \to H^0(K, V^*(1))^* \to 0.
  \end{gather*} 
  Here $H^1_{\mathrm{s}}(K_v, V) = H^1(K_v, V)/H^1_\f(K_v, V) = H^1_{\f}(K_v, V^*(1))^*$.
  
  We now count dimensions. We have $h^2(K_v, V) = h^1_{\mathrm{s}}(K_v, V)$ for $v \nmid p$, so the local terms for $v \nmid p$ cancel out; and $h^0(G_{K, S}, V) - h^1(G_{K, S}, V) + h^2(G_{K, S}, V) = -d_-(T)$ by Tate's global Euler characteristic formula. Finally, local Tate duality gives $h^2(K_v, V) = h^0(K_v, V^*(1))$. Collecting terms therefore gives
  \begin{align*}
   h^1_{\mathrm{rel}}(K, T) - h^0(K, T) 
   & = h^1_{\mathrm{str}}(K, T^*(1)) - h^0(K, T^*(1)) \\ &+ d_-(T) + \sum_{v \mid p} h^0(K_v, T^*(1)).
  \end{align*}
  Under our simplifying hypotheses, most of these terms are zero and the formula simplifies to
  \[ h^1_{\mathrm{rel}}(K, T) = d_-(T) + h^1_{\mathrm{str}}(K, T^*(1)).\]
  So $h^1_{\mathrm{rel}}(K, T) \ge d_-(T)$, with equality if and only if $h^1_{\mathrm{str}}(K, T^*(1)) = 0$.
 \end{proof}

\section{The case \texorpdfstring{$d_-(T) = 1$}{d-(T) = 1}}

 One can check that classes forming an Euler system always lie in $H^1_{\mathrm{rel}}$. Hence, if $d_-(T) = 1$ and $\bfc$ is an Euler system for $T$ with $\bfc_K$ non-torsion, then one has a rather precise picture of the cohomology of $T$, at least after inverting $p$; the space $H^1_{\mathrm{rel}}(K, T) \otimes L$ is one-dimensional, and $\bfc_K$ is an $L$-basis vector of this space.
 
 This situation, where $d_-(T) = 1$, may seem rather special, but it in fact covers several of the most familiar Euler systems:
 \begin{itemize}
  \item The Euler system of \emph{cyclotomic units}: here $K = \QQ$ and $T = \Zp(1)$.
  \item The Euler system of \emph{elliptic units}: here $K$ is imaginary quadratic and $T$ is again $\Zp(1)$.
  \item The Euler system of \emph{Beilinson--Kato elements}: here $K = \QQ$ and $T = T_f^*(1)$ where $T_f$ is the representation attached to a modular form of weight $\ge 2$, so that $T$ has rank 2 and $\sigma_\infty$ acts via a matrix conjugate to $\left(\begin{smallmatrix} -1 & 0 \\ 0 & 1 \end{smallmatrix}\right)$.
 \end{itemize}
 
 However, there are not many more examples beyond these. The problem is that in practice ``most'' representations $T$ have approximately the same number of $+1$ and $-1$ eigenvalues for complex conjugation; so $d_-$ is usually about $\tfrac{1}{2}[K: \QQ] \rank_{\cO}(T)$, which will be much larger than 1 unless $K$ and $T$ are both small.
 
 In practice, one is usually interested in Selmer groups with more sophisticated local conditions at $p$, rather than the (rather crude) strict and relaxed local conditions. The ``right'' local condition was defined by Bloch and Kato, using $p$-adic Hodge theory. We impose the assumption that $T$ is de Rham at the places above $p$ (which is automatically satisfied for all representations arising from geometry, by deep comparison theorems due to Faltings and Tsuji).
 
  \begin{definition}[{\cite[\S 3.7]{blochkato90}}]
   For $v \mid p$, define submodules $H^1_\f(K_v, T) \subseteq H^1_\mathrm{g}(K_v, T) \subseteq H^1(K_v, T)$ by
   \begin{align*}
    H^1_{\f}(K_v, T) &= \operatorname{ker}\left(H^1(K_v, T) \to H^1(K_v, T \otimes \mathbf{B}_{\mathrm{cris}})\right) \\
    H^1_{\mathrm{g}}(K_v, T) &= \operatorname{ker}\left(H^1(K_v, T) \to H^1(K_v, T \otimes \mathbf{B}_{\dR})\right)
   \end{align*}
   where $\mathbf{B}_{\mathrm{cris}}$ and $\mathbf{B}_{\dR}$ are Fontaine's $p$-adic period rings. Define the global Bloch--Kato Selmer group by
   \[ 
    H^1_{\f}(K, T) = \{ x \in H^1(K, T): \loc_v(x) \in H^1_{\f}(K_v, T) \text{ for all $v$} \}.
   \]
  \end{definition}
  
  From the theorems above, we see that if $d_-(T) = 1$ and an Euler system $\bfc$ exists for $T$ with $\bfc_K$ non-torsion, then $H^1_{\f}(K, T)$ has dimension either 1 or 0, depending on whether or not $\loc_v(\bfc_K) \in H^1_{\f}(K_v, T)$ for all primes $v \mid p$.
  
  If $V$ has all Hodge--Tate weights\footnote{Our conventions are that the Hodge--Tate weight of the cyclotomic character is $+1$. The Hodge--Tate weights of the representation $T_f$ attached to a weight $k$ modular form $f$ are 0 and $1-k$, so the representation $T_f^*(1)$ appearing in the Beilinson--Kato Euler system has weights 1 and $k$.} $\ge 1$ at some prime $v \mid p$, then $H^1_\mathrm{g}(K_v, T) = H^1(K_v, T)$ \cite[Lemma 6.5]{berger03}; if we suppose also that $H^0(K_v, T^*(1)) = 0$, as in Proposition \ref{prop:h1bound}, then we even have $H^1_\f(K_v, T) = H^1(K_v, T)$, so the condition $\loc_v(\bfc_K) \in H^1_{\f}(K_v, T)$ is automatically satisfied. For instance, this applies to the Euler system of Beilinson--Kato elements if the modular form $f$ has level coprime to $p$.
  
  On the other hand, if the Hodge--Tate weights are not all $\ge 1$ at $v$, one expects that $\loc_v(\bfc_K)$ should only be in $H^1_\f(K_v, T)$ if some ``unlikely coincidence'' occurs. For instance, in the setting of the Beilinson--Flach elements one can use a \emph{twisting} construction to produce an Euler system $\bfc'$ for $T = T_f^*$ (without the twist 1). It follows from Kato's explicit reciprocity law that this twisted Euler system has $\loc_p(\bfc'_{\QQ}) \in H^1_{\f}(\Qp, T_f^*)$ if and only if $L(f, 1) = 0$.
  
\section{Higher rank Euler systems}

 If $d_-(T) > 1$, what should one expect? Naively, one might guess that it would be \emph{easier} to build Euler systems in this context, since $H^1_{\mathrm{rel}}$ is forced to be large by Proposition \ref{prop:h1bound}. However, this doesn't seem to be the case: when $d_- $ is large it seems to be hard to construct elements. 
 
 An intuitive explanation of this comes from the following observation: systematic constructions of  elements in global cohomology groups only seem to work well when those groups are 1-dimensional, because otherwise the class ``doesn't know where to go'' within the space, and collapses to zero. (We shall call this \emph{Gross' trap}, since the observation was apparently first made by Dick Gross in the analogous setting of Heegner points on elliptic curves of analytic rank $> 1$.)
 
 One suggestion for resolving this problem, due to Perrin-Riou \cite{perrinriou98}, is that the ``correct'' object to associate to a general $T$ is not a collection of classes in $H^1(F, T)$, but rather classes in exterior powers of these modules. She defined a \emph{rank $r$ Euler system}, for $r \ge 1$, to be a collection of classes
 \[ \bfc_F \in \bigwedge^r_{\cO[\Delta_F]} H^1(F, T) \quad\text{for $K \subseteq_{\f} F \subset \cK$,} \]
 where $\Delta_F = \Gal(F / K)$, satisfying the Euler system norm relations ($\star$). 
 
 For simplicity, we shall state Perrin-Riou's conjecture under an auxiliary assumption: that either $d_-(T) = d_+(T)$, or every real place of $K$ remains real in $\cK$. (For example, we could take $K = \QQ$ and $\cK = \bigcup_m \QQ(\mu_m)^+$, where $\QQ(\mu_m)^+$ is the totally-real subfield of $\QQ(\mu_m)$.) This avoids complications with ranks varying between different complex-conjugation eigenspaces.
  
 \begin{conjecture}[Perrin-Riou]
  \label{conj:PR}
  For any global Galois representation $T$ arising in geometry, there exists an Euler system of rank $d_-(T)$ for $T$, satisfying a precise relation to the values of the $L$-function $L(T^*(1), \chi, s)$ for finite-order characters $\chi$ of $\Gal(\cK / K)$.
 \end{conjecture}
 
 \begin{remark}
  The exact relation to $L$-values is somewhat technical to state; see Perrin-Riou's monograph \cite{perrinriou95}, or the overview in \cite[Chapter 8]{rubin00}.
 \end{remark}
 
 This notion of higher-rank Euler systems has been extensively studied since, but it has proved to be rather thorny to work with, for two reasons. 
 
 Firstly, there are serious technical difficulties arising from the complicated algebra of wedge powers of modules over $\cO[\Delta_F]$. This makes it difficult to prove an analogue of Theorem \ref{thm:rubin} for Euler systems of rank $> 1$. Recent work of Burns--Sano \cite{burnssano16} strongly suggests that a better theory may be obtained by replacing the wedge power $\bigwedge^r_{\cO[\Delta_F]} H^1(F, T)$ with $\bigcap^r_{\cO[\Delta_F]}H^1(F, T)$, where $\bigcap^r$ denotes the ``exterior bi-dual'' functor, defined for modules $M$ over a ring $R$ by 
 \[ \bigcap^r_R M = \Hom_R\left( \bigwedge^r_R \Hom_R(M, R), R\right).\]
 
 However, another (possibly more serious) stumbling block is that there are very few interesting examples of rank $r$ Euler systems known for $r > 1$ (in particular, none which are known to be related to values of $L$-functions). In particular, it is \emph{not} expected that the Euler systems predicted by Perrin-Riou's conjecture should be constructed by building $r$ invidual elements in some canonical way, and then wedging them together (except in special cases, such as when $T$ is a direct sum of smaller representations); such an approach would fall into Gross' trap.
  
 \begin{remark}
  One exception to this gloomy outlook is provided by ongoing work of Nekov\'a\v{r} and Scholl (surveyed in \cite{nekovarscholl16}). Assuming a certain conjecture, the \emph{plectic conjecture}, their method gives a construction of Euler systems of rank $[F : \QQ]$ for certain Galois representations arising in the \'etale cohomology of Shimura varieties associated to reductive groups over totally real fields $F$. However, the plectic conjecture is currently wide open. 
  
  Another, unrelated approach is due to Urban, who has devised a method of constructing higher-rank Euler systems via Eisenstein congruences; but this approach (as presently formulated) requires one to assume bounds on congruence ideals as input to the method, and these congruence ideals are closely related to Selmer groups, so using these classes as input to a version of Theorem \ref{thm:rubin} would result in a circular argument.
 \end{remark}

\section{Euler systems with local conditions}

 In 2014, in joint work with Lei, we discovered a new example of an Euler system:
 
 \begin{theorem}[{\cite[Corollary 6.4.5]{leiloefflerzerbes14}}]
  \label{thm:LLZ}
  Let $f, g$ be two modular forms of weight $2$ and prime-to-$p$ level, and let
  \[ T = (T_f \otimes T_g)^*, \]
  where $T_f$ and $T_g$ are the Galois representations attached to $f$ and $g$. Then there exists a collection of classes $\bfc_{\QQ(\mu_m)} \in H^1(\QQ(\mu_m), T)$ satisfying compatibility relations close to $(\star)$.
 \end{theorem}
 
 \begin{remark}
  This theorem is, of course, vacuous as stated, since the $\bfc_{\QQ(\mu_m)}$ could all be 0; but we can also show in many cases that $\bfc_\QQ$ is non-torsion. 
 \end{remark}
 
 There are several curious features of the Euler system of Bei\-lin\-son--Flach elements. Firstly, it has the ``wrong'' rank:  $T$ is 4-dimensional and odd, so $d_-(T) = 2$. Thus Conjecture \ref{conj:PR} would predict a rank 2 Euler system, not a rank 1 Euler system.
 
 Secondly, the norm-compatibility relations satsified by the Bei\-lin\-son--Flach elements for $\ell = p$ are not the expected ones. If we write $\QQ_r = \QQ(\mu_{p^r})$, then we obtain formulae of the form
 \[ \operatorname{cores}_{\QQ_r}^{\QQ_{r+1}}(\bfc_{\QQ_{r + 1}}) = (\alpha_f \alpha_g) \cdot \bfc_{\QQ_r}, \]
 where $\alpha_f$ and $\alpha_g$ are some choices of roots of the Hecke polynomials of $f$ and $g$ at $p$. If $f$ and $g$ are ordinary, we may choose $\alpha_f$ and $\alpha_g$ to be $p$-adic units; then we can re-normalise by setting $c_{\QQ_r}' = (\alpha_f \alpha_g)^{-r} \bfc_{\QQ_r}$ to obtain the expected Euler system relation. However, if $\alpha_f$ and $\alpha_g$ are not $p$-adic units, then there is no way to re-normalise the elements $\bfc_{\QQ_r}$ to be norm-compatible without introducing denominators. 
 
 It turns out that these distorted norm-compatibility relations at $p$ are \emph{unavoidable}. The Beilinson--Flach classes are automatically in $H^1_{\mathrm{g}}$, since they are constructed geometrically; and $T$ has Hodge--Tate weights $\{0, 1, 1, 2\}$, which are \emph{not} all $\ge 1$. This means there is a local obstruction to having norm-compatible classes, because of the following theorem of Berger:
 
 \begin{theorem}[{\cite[Theorem A]{berger05}}]
  Let $T$ be an irreducible $\cO$-linear de Rham representation of $G_{K_v}$ of dimension $> 1$, for $K_v$ a $p$-adic field, and suppose we are given classes $x_n \in H^1_\mathrm{g}(K_v(\mu_{p^n}), T)$ for all $n \ge 1$ which are compatible under corestriction. 
  
  Then either $T$ has all Hodge--Tate weights $\ge 1$, or $x_n = 0$ for all $n$.
 \end{theorem}
 
  So if $T|_{G_\Qp}$ is irreducible (which can occur) then any collection of norm-compatible classes lying in $H^1_{\f}$ at $p$ would either have to localise to 0 at $p$ (which is unlikely, because we expect the strict Selmer group to be generically 0); or it would have to have a denominator growing in the cyclotomic tower, at a certain minimum rate determined by the valuation of $\alpha_f \alpha_g$. This is exactly the behaviour one sees for the Beilinson--Flach classes.
 
 Fortunately, for the machinery of Kolyvagin derivatives, one is mainly interested in classes over $\QQ(\mu_m)$ where $m$ is a squarefree product of primes coprime to $p$, so this ``distortion'' of the $p$-direction norm relations does not rule out applications to Selmer groups.  One can use this to show (under the usual auxillary ``big image'' hypotheses) that when $\bfc_\QQ$ is non-torsion, the group $H^1_{\mathrm{str}}(\QQ, T^*(1))$ is finite, and $H^1_{\f}(\QQ, T)$ is of rank 1 and is spanned by $\bfc_{\QQ}$ after inverting $p$.
 
 \begin{remark}
  In the three classical examples of Euler systems listed in the previous section, the cohomology classes are also constructed geometrically, so they likewise lie in $H^1_{\mathrm{g}}$; but in these examples the Hodge--Tate weights are all $\ge 1$, so $H^1_{\mathrm{g}}$ is the whole of the local cohomology at $p$ and Berger's theorem is no obstruction. The novel feature of the Beilinson--Flach classes is that they are in $H^1_{\mathrm{g}}$ at $p$ in a situation where this is a nontrivial condition.
 \end{remark}
  
\section{A conjecture}

 These properties of the Beilinson--Flach elements suggests that Perrin-Riou's conjecture \ref{conj:PR} is not the whole story. This motivates a more general Euler system conjecture, which we explain below.
 
 For technical reasons, the conjecture is simplest to state if we abandon the assumption that the coefficient field $L$ is a finite extension of $\Qp$, and instead assume that it is a finite extension of $\operatorname{Frac} W(\overline{\FF}_p)$, where $W(-)$ denotes Witt vectors. (This base-extension is not needed if $K = \QQ$.) As before, we write $\cO$ for the ring of integers of $L$.
 
 Let $K$,  $\cK$ and $T$ be as in \S \ref{sect:ES} above, and write $V = T \otimes_\cO L$. We assume $V$ is unramified almost everywhere and de Rham at the places above $p$. We also assume that $V$ is irreducible and that $H^0(K, V) = H^0(K, V^*(1)) = 0$.
 
 \begin{definition}
  We define
  \begin{align*}
   r_0(T) &\coloneqq d_-(T) - \sum_{v \mid p} \dim_L \Fil^0 \DD_{\dR}(K_v, V),\\ 
   r(T) &\coloneqq  \max\left(0, r_0(T)\right).
  \end{align*}
 \end{definition}
 
 One checks easily\footnote{One has $d_-(T) + d_-(T^*(1)) = d_-(T) + d_+(T) = [K:\QQ] \dim V$; while for each $v \mid p$, there is a perfect pairing $\DD_{\dR}(K_v, V) \times \DD_{\dR}(K_v, V^*(1)) \to L$, and the two $\Fil^0$'s are orthogonal complements, so their dimensions sum to $[K_v : \Qp] \dim V$.} that $r_0(T^*(1)) = -r_0(T)$, so for any $T$, at least one  of $r(T)$ and $r(T^*(1))$ is zero. An application of Poitou--Tate duality gives the following relation:
 
 \begin{proposition}
  We have
  \[ h^1_{\f}(K, T) - r(T) = h^1_{\f}(K, T^*(1)) - r(T^*(1)).\]
  In particular, if $r(T^*(1)) = 0$, then we have 
  \[ h^1_{\f}(K, T) \ge r(T),\]
  with equality if and only if $h^1_\f(K, T^*(1)) = 0$.\qed
 \end{proposition}
 
 The significance of $r(T)$ is as follows. The Bloch--Kato conjecture predicts that we should have
 \[ h^1_{\f}(K, T) = \ord_{s = 0} L(T^*(1), s).\]
 On the other hand, Deligne has defined an archimedean $L$-factor $L_\infty(T^*(1), s)$, which is a product of $\Gamma$-functions depending on the Hodge--Tate weights of $T$ and the action of complex conjugation on it. Deligne's conjectures predict that (under our hypotheses on $T$) the function
 \[ \Lambda(T^*(1), s) \coloneqq L(T^*(1), s) L_\infty(T^*(1), s) \]
 should be meromorphic on $\CC$, and holomorphic at $s = 0$. The archimedean factor $L_\infty(T^*(1), s)$ has no zeroes, but it does have poles, and $r(T)$ is exactly the order of the pole of $L_\infty(T^*(1), s)$ at $s = 0$. Hence, if $\Lambda(T^*(1), s)$ is to be holomorphic at $s = 0$, the function $L(T^*(1), s)$ must vanish there to order at least $r(T)$. In other words, $r(T)$ is the ``Archimedean contribution'' to the order of vanishing of $L(T^*(1), s)$. 
   
 \begin{remark}
  It is expected that for almost all values of $s$ (whenever $T$ does not have ``motivic weight $-1$'') the functional equation will force $\Lambda(T^*(1), s)$ to be non-vanishing at $s = 0$; so this Archimedean contribution should actually completely determine the order of vanishing.
 \end{remark}
 
 \begin{definition}
  For an integer $r \ge 0$, we say $T$ is \emph{$r$-critical} if $r(T) = r$ and $r(T^*(1)) = 0$.
 \end{definition}
 
 The second condition is, of course, redundant if $r > 0$; it is included only in order to ensure that $0$-critical agrees with the usual notion of critical, which is that neither $L_\infty(T, s)$ nor $L_\infty(T^*(1), s)$ has a pole at $s = 0$.
 
 We can now formulate our first conjecture on the existence of Euler systems. We first consider only fields unramified above $p$, postponing discussion of the ``$p$-direction'' until later.
 
  \begin{conjecture}[rough form]
   \label{conj:basicconjecture}
   If $T$ is $r$-critical, there exists a collection of cohomology classes
   \[ \bfc_F \in \bigwedge^r H^1_{\f}(F, T),\]
   where $F$ varies over finite extensions of $K$ inside $\cK$ that are unramified above $p$, satisfying the Euler system compatibility relation ($\star$); and the bottom class $\bfc_K$ is non-zero if and only if $L^{(r)}(T^*(1), 0) \ne 0$. 
  \end{conjecture}
  
  \begin{remark}
   This conjecture is not precise, since we have not attempted to formulate a relation to $L$-functions. This should be roughly as follows: suppose $T$ is the $p$-adic realisation of a motive. Then Beilinson's conjecture predicts that $L^{(r)}(T^*(1), 0)$ should be given by Beilinson's regulator map applied to an element in the $r$-th wedge power of a motivic cohomology group, and it is natural to expect that $c_K$ should be the $p$-adic realisation of this motivic element.
   
   It is also very possible that the conjecture may need some modification to account for denominators, replacing $\bigwedge^r H^1_{\f}(F, T)$ with some larger lattice in $\bigwedge^r H^1_{\f}(F, V)$, such as the exterior bi-dual lattice $\bigcap^r H^1_{\f}(F, T)$, as in the work of Burns--Sano cited above. However, we shall not pursue this here, since we want to focus primarily on cases where $r(T) = 1$; in this case the ``naturally occurring'' elements do indeed seem to lie in $H^1_{\f}(F, T)$.
  \end{remark}
  
  For instance, if $f, g$ are weight 2 modular forms, then we have
  \[ r((T_f \otimes T_g)^*(m)) = 
   \begin{cases}
    2 & \text{if $m \ge 1$}\\
    1 & \text{if $m = 0$}\\
    0 & \text{if $m \le -1$}.
   \end{cases}
  \]
  Thus our conjecture predicts that there should be Euler systems of \emph{multiple} ranks attached to different twists of $T = (T_f \otimes T_g)^*$. There should be Euler systems of rank 2 attached to $T(m)$ for each $m \ge 1$, which are the objects predicted by Perrin-Riou's conjecture; but there should also be a rank 1 Euler system for $T$ itself, which is the Euler system of Beilinson--Flach elements. We shall consider this example in more detail below. 
  
  It is important to note that this conjecture is not, in itself, particularly novel; for instance, one can deduce it from Perrin-Riou's Conjecture \ref{conj:PR}, by applying various linear functionals to the conjectural rank $d_-$ Euler system to move it down to rank 1, as we shall describe in a later section. 
  The reason why we feel that Conjecture \ref{conj:basicconjecture} is interesting is that it may be more approachable than Conjecture \ref{conj:PR}. We optimistically hope that when our conjecture predicts a \emph{rank 1} Euler system (i.e.~when we have a geometric Galois representation with $r(T) = 1$) then one can reasonably expect to construct the necessary cohomology classes directly.
  
  Moreover, the lower-rank Euler systems predicted by Conjecture \ref{conj:basicconjecture} still have powerful arithmetic applications. Although they have lower ranks than those predicted by Perrin-Riou, this is ``compensated for'' by their additional local property at $p$ -- namely, they lie in $H^1_{\f}$. As shown in \cite[Appendix B]{leiloefflerzerbes14b}, when $r(T) = 1$ one can adapt the proof of Theorem \ref{thm:rubin} to make use of this additional information:
  
  \begin{proposition}
   Suppose $r(T) = 1$ and there exists a rank 1 Euler system for $T$ such that $\bfc_F \in H^1_{\f}(F, T)$ for all $F$ and $\bfc_K$ is non-torsion. Under some auxilliary technical hypotheses, then $H^1_{\f}(K, T^*(1))$ is finite, $H^1_{\f}(K, T)$ has rank 1 and is spanned by $c_K$, and $H^1(K, T)$ has rank $d_-(V)$.
  \end{proposition}
   
  The case $r = 0$ of Conjecture \ref{conj:basicconjecture} is not at all trivial. It predicts the existence of collections of elements of the group rings $\cO[\Delta_F]$ satisfying some norm-compatibility properties; and the expected relation to $L$-values simplifies greatly in this case, predicting that the image of the element $c_F \in \cO[\Delta_F]$ under evaluation at a character $\chi$ of $\Delta_F$ should give the critical $L$-value $L(T^*(1), \chi, 0)$ divided by an appropriate period. 
  
  There are several naturally-occurring examples of such elements: for instance, one has the Stickelberger elements attached to $T = \cO(\chi)$, where $\chi$ is a Dirichlet character with $\chi(-1) = 1$, and the Mazur--Tate elements for $T = T_f(1)$ where $f$ is a weight 2 modular form.
  
 \section{Ordinarity conditions at $p$}
 
  We now consider the question of norm relations in the $p$-direction. If $r(T) < d_-(T)$, so that our conjecture predicts Euler systems of ``non-optimal'' rank, then there must be at least one prime above $p$ at which $V$ has a Hodge--Tate weight $\le 0$. So Berger's theorem shows that there is an obstruction to having norm-compatible systems of geometric classes over the $p$-cyclotomic tower. In other words, we should not expect to have such an interpolation unless the local representations are reducible. 
  
  In fact, it turns out that we need subrepresentations of a very specific kind:
  
  \begin{definition}
   Let $v$ be a prime above $p$. A \emph{Panchishkin subrepresentation} of $V$ at $v$ is a subspace $V_v^+ \subseteq V$ such that
   \begin{itemize}
    \item $V_v^+$ is stable under $G_{K_v}$,
    \item $V_v^+$ has all Hodge--Tate weights $\ge 1$,
    \item $V / V_v^+$ has all weights $\le 0$.
   \end{itemize}
  \end{definition}
 
  Note that $V_v^+$ is unique if it exists. If such a $V_v^+$ exists, then (up to minor grains of salt), one sees that $H^1_{\f}(K_v, V)$ is simply the image of the natural map $H^1(K_v, V_v^+) \to H^1(K_v, V)$.
  
  \begin{definition}
   We say $V$ satisfies the \textbf{rank $r$ Panchishkin condition} if $r(V) = r$, $r(V^*(1)) = 0$, and Panchishkin subrepresentations $V_v^+$ exist for all $v \mid p$. 
  \end{definition}
  
  Note that if this holds, we must necessarily have $\sum_{v \mid p} [K_v: \Qp] \dim_L(V_v^+) = d_+(V) + r$.
  
  This condition was introduced in the case $r = 0$ by Panchishkin, who suggested that the rank 0 Panchishkin condition was the ``correct'' condition for a (bounded) $p$-adic $L$-function to exist -- in other words, for rank 0 Euler systems to interpolate in the $p$-cyclotomic tower. 
  
  \begin{remark}
   The Panchishkin condition is closely related to the notion of \emph{ordinarity}. This has various formulations, but one flavour is to require that $V |_{G_{K_v}}$ have a decreasing filtration by subrepresentations $V_v^{(i)}$ such that each quotient $V_v^{(i)} / V_v^{(i+1)}$ has all Hodge--Tate weights equal to $i$. Thus $V$ is ordinary at some prime $v \mid p$ if and only if all its Tate twists $V(j)$ have Panchishkin subrepresentations. However, full ordinarity of this kind is a rather restrictive condition, and (as we shall see later) it is interesting and instructive to see how much of this condition is actually relevant in specific situations.
  \end{remark}
 
  \begin{conjecture}
   \label{conj:pdirection}
   If $T$ is $r$-critical and satisfies the rank $r$ Panchishkin condition, then there should be a collection of classes $\bfc_F \in \bigwedge^r H^1(F, T)$ as in Conjecture \ref{conj:basicconjecture} for all $K \subseteq_{\f} F \subset \cK$ (not just those unramified above $p$).
  \end{conjecture}
 
  Notice that if $r = d_-(V)$, then the rank $r$ Panchishkin condition is trivially satisfied (since we can take $V_v^+ = V$ for every $v \mid p$). This is why ordinarity plays no role in the Euler system of Kato, for instance; but for Euler systems of non-optimal rank, the Panchishkin condition is a non-trivial restriction.
  
  \subsection{Example A: Rankin--Selberg convolutions}
   \label{sect:exRS}
   
   Consider the representations $T = (T_f \otimes T_g)^*(m)$ introduced in Theorem \ref{thm:LLZ}, for $f, g$ modular forms of weights $k+2,\ell + 2$, with $k,\ell \ge 0$. Note that $d_-(T) = 2$. We assume $k \ge \ell$ without loss of generality.
   
   \begin{itemize}
    
    \item When $m \ge 1$, the representation $T$ is $2$-critical; so Conjecture \ref{conj:basicconjecture} predicts a rank 2 Euler system, and the Panchishkin condition is automatic, so this Euler system should extend up the $p$-cyclotomic tower without further hypotheses.
    
    \item When $0 \ge m \ge -\ell$, the representation is 1-critical; in this case, we need to take $V_v^+$ to be a 3-dimensional subrepresentation of $V |_{G_{\Qp}}$, i.e.~the orthogonal complement of a 1-dimensional subrepresentation of $V_f \otimes V_g$ of Hodge--Tate weight 0. 
    
    If we assume $f$ and $g$ are both ordinary, then $V_f$ and $V_g$ both have one-dimensional subrepresentations $V_f^+$ and $V_g^+$ (each of which is unramified, with Hodge--Tate weight 0) and we can take the 1-dimensional sub to be $V_f^+ \otimes V_g^+$.
    
    \item When $-1-\ell \ge m \ge -k$, the representation is 0-critical. Hence, in order to find a rank 0 Euler system in the $p$-direction -- that is, a $p$-adic $L$-function -- we require the existence of a 2-dimensional subrepresentation of $V_f \otimes V_g$ accounting for the two highest Hodge--Tate weights $\{0, -1-\ell\}$. Such a subrepresentation exists when $f$ has strictly larger weight, i.e.~$k > \ell$, and $f$ is ordinary: we can take $V_f^+ \otimes V_g$. Note that we do not need to assume any ordinarity condition on $g$ here.
   
   \end{itemize}
  
   (We do not need to consider $m \le -1-k$, since then $r(T^*(1))$ is no longer zero and our conjecture does not apply.)
   
   So we should expect a rank 1 Euler system in the $p$-direction when both $f$ and $g$ are ordinary; but to form a $p$-adic $L$-function, we only need to assume ordinarity for whichever of the two forms has the highest weight. This matches exactly the behaviour one observes for Beilinson--Flach elements and the Panchishkin--Hida $p$-adic Rankin--Selberg $L$-function.
  
  \subsection{Example B: The spin representation for \texorpdfstring{$\operatorname{GSp}(4)$}{GSp4}}

   We now consider a more sophisticated example. We take $\cF$ a cuspidal Siegel modular eigenform of genus 2 and weight 3. By work of Taylor and Weissauer \cite{weissauer05}, this gives rise to a Galois representation
   \[ \rho_{\cF}: G_{\QQ} \to \operatorname{GSp}_4(\overline{\QQ}_p).\]
   Composing this with the canonical inclusion of $\operatorname{GSp}_4$ into $\GL_4$ gives a 4-dimensional representation of $G_{\QQ}$, which we denote by $V_{\cF}$. 
   
   \begin{remark}
    This representation is called the \emph{spin representation}, for reasons which only become obvious when one considers more general symplectic groups $\operatorname{GSp}_{2g}$. The Langlands dual of $\operatorname{GSp}_{2g}$ is the \emph{spin similitude group} $\operatorname{GSpin}_{2g+1}$, which acts naturally on a $2^g$-dimensional space of ``spinors''. However, for $g = 2$ there is an exceptional isomorphism $\operatorname{GSpin}_{5} \cong \operatorname{GSp}_{4}$, and the spinor space is simply the 4-dimensional defining representation of $\operatorname{GSp}_4$. 
    
    The spin Galois representation should not be confused with the \emph{standard} representation, given by composing $\rho_{\cF}$ with the 5-dimen\-sion\-al defining representation of $\operatorname{SO}_5 \cong \operatorname{PGSp}_4$.
   \end{remark} 
   
   If $p$ does not divide the level of $\cF$, the local behaviour of $\rho_{\cF}$ at $p$ is determined by the Hecke eigenvalues of $\cF$ at $p$. The Hecke algebra has two generators, corresponding to the double cosets
   \[ T(p) = \left[ \left( \begin{smallmatrix} 1 \\& 1 \\&& p \\&&&p  \end{smallmatrix}\right)\right] \quad\text{and}\quad 
   T_1(p^2) = \left[ \left( \begin{smallmatrix} 1 \\& p \\&& p \\&&&p^2  \end{smallmatrix}\right)\right].\]
   
   These correspond, respectively, to the two maximal proper parabolic subgroups of $\operatorname{GSp}_4$: the Siegel and Klingen parabolics. We say $\cF$ is \emph{Siegel-ordinary} if $T(p)$ acts as a $p$-adic unit, and \emph{Klingen-ordinary} if $T_1(p^2)$ does so.
   
   \begin{theorem}[{\cite[Corollary 1]{urban05}}] \
    \begin{enumerate}
     \item[(i)] If $\cF$ is Siegel-ordinary, then $\rho_{\cF}(G_{\Qp})$ stabilises a line in $V_{\cF}$.
     \item[(ii)] If $\cF$ is Klingen-ordinary, then $\rho_{\cF}(G_{\Qp})$ stabilises a plane in $V_{\cF}$.
    \end{enumerate}
   \end{theorem}
   
   \begin{remark}
    Urban proves (ii) under an additional technical condition, that the automorphic representation $\Pi$ generated by $\cF$ be ``stable at $\infty$'' (see Remark (i) \emph{loc.cit.}). This hypothesis can now be removed, as a consequence of Arthur's classification of cuspidal automorphic representations of $\operatorname{GSp}_4$, announced in \cite{arthur04} and proved in \cite{geetaibi18}.
   \end{remark}
      
   What does our conjecture say in this case? The representation $V_{\cF}$ has $d_-(V) = 2$, and Hodge--Tate weights $\{0, -1, -2, -3\}$. Setting $V = V_{\cF}^*(-j)$, we expect that:
   
   \begin{itemize}
    \item when $j \le -1$, $V$ is $2$-critical, and we expect a rank 2 Euler system in the $p$-direction without any ordinarity conditions;
    
    \item when $j = 0$, $V$ is $1$-critical, so we expect a rank 1 Euler system, and if we wish to extend this in the $p$-direction, we need to assume $\cF$ is Siegel-ordinary; 
     
    \item when $j = 1$, $V$ is $0$-critical (i.e.~critical in the sense of Deligne), so we expect a rank 0 Euler system; and the condition required to interpolate this into a $p$-adic $L$-function is that $\cF$ should be Klingen-ordinary.
   \end{itemize}
  
   More generally, this analysis goes over to Siegel modular forms of any cohomological weight, and one again finds that Siegel-ordinarity is the condition for a rank 1 Euler system, and Klingen-ordinarity the right condition for a $p$-adic $L$-function. This is exactly what one sees in two recent papers: our work with Skinner on the construction of a (rank 1) Euler system for these representations \cite{loefflerskinnerzerbes17}; and work of Dimitrov, Januszweski and Raghuram on the construction of a $p$-adic $L$-function \cite{DJR18}.
 
 \section{Iwasawa theory and Greenberg Selmer groups}
  
  Let $F_\infty = \bigcup_{n \ge 1} F_n$ be a $\ZZ_p^m$-extension of $K$ contained in $\cK$, for some $m \ge 1$. We assume $F_\infty$ contains the cyclotomic $\Zp$-extension $K_\infty / K$. Let $\Gamma = \Gal(F_\infty / K)$, and let $\Lambda(\Gamma)$ be the Iwasawa algebra of $\Gamma$ with coefficients in $\cO$. 
  
  \subsection{The rank 0 case}
  
   For representations $V$ satisfying the rank 0 Panchishkin condition, we expect that there should be a $p$-adic $L$-function, which should be an element of $\Lambda(\Gamma)$ interpolating $L$-values $L(V^*(1), \chi, 0) / \text{(period)}$ as $\chi$ varies over finite-order characters of $\Gamma$. 
   
   In a ground-breaking paper \cite{greenberg89} in Iwasawa's 70th birthday proceedings, Greenberg showed how define a Selmer group associated to $V$ over $F_\infty$, and thus formulate an Iwasawa main conjecture. He introduced the following two objects, known as ``Greenberg Selmer groups'':
   
   \begin{itemize}
    \item a subgroup $H^1_{\operatorname{Gr}}(F_\infty, T)$ of $H^1_{\Iw}(F_\infty, T) = \varprojlim_n H^1(F_n, T)$ defined by local conditions, in which the local condition at $v \mid p$ is the image of the cohomology $H^1_{\Iw}(F_{\infty, v}, T \cap V_v^+)$ where $V_v^+$ is the Panchishkin subrepresentation;
    
    \item a subgroup $H^1_{\operatorname{Gr}}(F_\infty, T^\vee(1))$ of $H^1(F_\infty, T^\vee(1))$, where $T^\vee = \Hom(T, \Qp/\Zp)$, defined similarly using the orthogonal complement of $T \cap V_v^+$ in $T^\vee$.
   \end{itemize}
  
   The compact Greenberg Selmer group $H^1_{\operatorname{Gr}}(F_\infty, T)$ is a finitely-generated $\Lambda(\Gamma)$-module, and the discrete version $H^1_{\operatorname{Gr}}(F_\infty, T^\vee(1))$ is a co-finitely-generated one (i.e. its Pontryagin dual $\mathfrak{X}(F_\infty, T)$ is finitely generated). Moreover, the ranks of $H^1_{\operatorname{Gr}}(F_\infty, T)$ and $\mathfrak{X}(F_\infty, T)$ are the same, by a Poitou--Tate duality computation. Greenberg's main conjecture is that these modules are both torsion, and that the characteristic ideal of $\mathfrak{X}(F_\infty, T)$ is generated by the $p$-adic $L$-function.
  
  \subsection{Higher ranks}
   \label{sect:mainconj}
  
   How should this look for $r$-critical representations, when $r > 0$? If the rank $r$ Panchishkin condition holds, the definitions of the two Greenberg Selmer groups still make sense; but one finds that their ranks differ by $r$. Moreover, if conjecture \ref{conj:pdirection} holds, then the Euler system classes $\bfc_{F_n}$ for $n \ge 1$ define an element $\bfc_{F_\infty}$ of $\bigwedge^r H^1_{\operatorname{Gr}}(F_\infty, T)$. The natural conjecture appears to be that the quotient
   \[ \frac{\bigwedge^r H^1_{\operatorname{Gr}}(F_\infty, T)}{\Lambda(\Gamma) \cdot \bfc_{F_\infty}} \]
   should be torsion as a $\Lambda$-module, and that its characteristic ideal should coincide with that of $\mathfrak{X}(F_\infty, T)$. When $r = d_-(T)$, so that the local conditions in the Greenberg Selmer groups are the trivial ones, this conjecture has already been formulated by Perrin-Riou; see chapter 8 of \cite{rubin00}. However, as explained above, we feel that settings with $r = 1$ may be more approachable.
   
 \section{Rank-lowering operators and reciprocity laws}
 
  There exist ``twisting'' operators for Euler systems (of any rank): if $\chi$ is a continuous character of $\Gal(\cK / K)$ unramified outside $\cN$, then there is a canonical bijection between Euler systems for $T$ and for $T(\chi)$. See e.g.~\cite[\S 6.3]{rubin00}. In particular, if $\cK$ contains the $p$-power cyclotomic extension $K(\mu_{p^\infty})$, then an Euler system for $T$ is also an Euler system for all of its Tate twists $T(n)$.
  
  How do these twisting maps interact with the predictions of Conjecture \ref{conj:pdirection}? Let us suppose that $T$ is $r$-critical, and $T(\chi)$ is $s$-critical for some integers $r \ge s$; we would like to compare the conjectured Euler systems for $T$ and for $T(\chi)$. Let us write $T_v^+$ for the Panchishkin subrepresentations for $T$, and $T_v^{++}$ for those\footnote{More precisely, $T_v^{++}$ is the subrepresentation of $T$ such that $T_v^{++}(\chi)$ is the Panchishkin subrepresentation of $T(\chi)$ at $v$.} of $T(\chi)$. 
  
  Our assumptions imply that
  \[ \sum_{v \mid p} [K_v: \Qp] \rank(T_v^{++}) \le \sum_{v \mid p}[K_v: \Qp] \rank(T_v^{+}),\]
  and it seems reasonable to expect a relation whenever $T_v^{++} \subseteq T_v^+$ for all $v \mid p$.
  
  \subsection{The equal-rank case}
  
   If $r = s$, then this condition will force $T_v^{+} = T_v^{++}$ for all $v$; and one can reasonably expect that the rank $r$ Euler systems associated to these two representations by \ref{conj:pdirection} should coincide under twisting. This gives the following refinement of Conjecture \ref{conj:pdirection}:
   
   \begin{conjecture}
    \label{conj:prefinement}
    Suppose we are given a collection $\cP$ of local subrepresentations $T_v^{+} \subseteq T |_{G_{K_v}}$ for all $v \mid p$, with $r(\cP) = -d_+(T) +  \sum_{v} [K_v: \Qp] \rank T_v^+ \ge 0$. Let $\Sigma(\cP)$ be the set of characters $\chi: \Gal(\cK / K) \to \overline{L}^\times$, unramified outside $\cN$ and de Rham above $p$, such that $T(\chi)$ is $r$-critical and $T_v^{+}(\chi)$ is a Panchishkin subrepresentation of $T(\chi) |_{G_{K_v}}$ for all $v \mid p$.

    If $\Sigma(\cP)$ is non-empty, there exists an Euler system $\bfc(\cP)$ for $T$ of rank $r = r(\cP)$ such that for every $\chi \in \Sigma(\cP)$ and every field $F$ with $\cK \supset F \supseteq_{\f} K$, the image of $\bfc(\cP)$ in $\bigwedge^r H^1(F, T(\chi))$ is the class $\bfc_F$ predicted by Conjecture \ref{conj:pdirection} applied to $T(\chi)$.
   \end{conjecture}
  
   In other words, the Euler system depends not on the specific twist $\chi$ that we choose, but only on which local subrepresentations are the Panchishkin subrepresentations for $\chi$.
   
   For $r = 0$, this is a familiar property of $p$-adic $L$-functions -- that a single $p$-adic $L$-function will often interpolate critical values of twists with a range of infinity-types, as long as these twists are all critical ``in the same way'', i.e. they all admit the same Panchishkin subrepresentation. For $r = 1$, the refined conjecture implies compatibilities under twisting between cohomology classes arising from very different geometric constructions. In the case of the Beilinson--Flach elements, this compatibility does indeed hold \cite[Theorem 6.3.4]{KLZ17}, but it is far from easy to show; it seems to be a rather deep result, requiring the full force of Kings' theory of $p$-adic interpolation of Eisenstein classes.
    
  \subsection{Rank-lowering}
  
   We now suppose that $T_v^{++} \subseteq T_v^+$ for all $v$ and that we have strict inequality $r > s$. Let $t = r-s$, and for each $v$ set $T_v^\sharp = T_v^+ / T_v^{++}$. We have $\sum_v [K_v : \Qp] \rank(T_v^\sharp) = t > 0$, so at least one of the $T_v^\sharp$ is non-zero. Let $F_\infty$ be a $p$-adic Lie extension of $K$ inside $\cK$, chosen such that $F_\infty$ contains the cyclotomic $\Zp$-extension $K_\infty$, and $\chi$ factors through $\Gal(F_\infty / K)$. Write $H^1_+(F_\infty, T)$ for the kernel of the map $H^1_\Iw(F_\infty, T) \to \bigoplus_{v \mid p} H^1_{\Iw}(F_{\infty, v}, T / T_v^+)$, and similary $H^1_{++}(F_\infty, T)$. Then there is an exact sequence
  \[ 
   0 \to H^1_{++}(F_\infty, T) \to H^1_+(F_\infty, T) \to \bigoplus_{v \mid p} H^1_{\Iw}(F_{\infty, v}, T_v^\sharp). 
  \]
  The final group in this sequence, however, is rather simpler than the previous ones, since it depends only on local information at $p$; in particular its rank over the Iwasawa algebra $\Lambda$ is known -- it is exactly $t$. Moreover, the \emph{local epsilon-isomorphism conjecture} of Fukaya and Kato \cite[Conjecture 3.4.3]{fukayakato06} predicts that its top wedge power should be canonically identified with $I \otimes \bigwedge^{t}_{\cO}(\bigoplus_v T_v^\sharp)$, where $I$ is a certain explicit fractional ideal in $\Lambda(\Gamma)$ (which is the unit ideal unless one of the local $L$-factors associated to the $T_v^\sharp$ has an exceptional zero). Note that our assumption that $\cO$ contain $W(\overline{\FF}_p)$ is essential here.
  
  Sadly, for general $T_v^\sharp$ and $F_\infty$ this local conjecture appears to be out of reach; but it is known in the important special case when $K_v$ is unramified over $\Qp$, the local extension $F_{\infty, v}$ is abelian over $\Qp$, and $T_v^\sharp$ is crystalline (see \cite{benoisberger08} or \cite{LVZ}). In this case, the required trivialisation is given by the determinant of Perrin-Riou's regulator map
  \[ \mathcal{L}_V : H^1_{\Iw}(F_{\infty, v}, T^\sharp) \to \mathcal{H}(\Gamma) \otimes \DD_{\mathrm{cris}}(K_v, V), \]
  where $\mathcal{H}(\Gamma)$ is the algebra of locally-analytic distributions on $\Gamma = \Gal(F_{\infty, v} / K_v)$.
  
  In cases when we can establish the local $\varepsilon$-isomorphism conjecture, we obtain a supply of linear functionals $\bigwedge^t H^1_{+}(F_{\infty, v}, T) \to \Lambda$. Recalling that $t = r - s$, these functionals can be regarded as linear maps
  \[ \bigwedge^r H^1_{+}(F_{\infty, v}, T) \to \bigwedge^s H^1_{+}(F_{\infty, v}, T), \]
  whose image is actually contained in $\bigwedge^s H^1_{++}(F_{\infty, v}, T)$. 
  
  Allowing the field $F_\infty$ to vary over $p$-adic Lie extensions of $F$ inside $\cK$, we obtain a map from Euler systems of rank $r$ for $T$ with local conditions given by $T_v^+$, to Euler systems of rank $s$ for $T(\chi)$ with local conditions given by $T_v^{++}$. We can now make the (rather optimistic) conjecture that the Euler systems predicted by Conjecture \ref{conj:pdirection} should be compatible under these ``rank-lowering operators''.
  
  \subsection*{The case $s = 0$} Let us now home in on the case $s = 0$ for a moment. We have already noted that our rank 0 Euler systems for $T(\chi)$ should be families of elements of group rings $\cO[\Delta_F]$, interpolating the critical values $L(T^*(1)(\chi^{-1}), \tau, 0)$ as $\tau$ varies over finite-order characters of $\Delta_F$. Compatible systems of such objects, as $F$ varies over subfields of $F_\infty$, can thus be regarded as $p$-adic $L$-functions. So our ``rank-lowering'' conjecture predicts that a map from rank $r$ Euler systems to rank 0 Euler systems, given (essentially) by the $r$-th wedge power of the Perrin-Riou regulator map, should send the rank $r$ Euler systems predicted by Conjecture \ref{conj:pdirection} to $p$-adic $L$-functions interpolating the critical values of twists of $T$. 
  
  Results of this kind -- relating Euler systems to critical $L$-values -- are generally known as ``explicit reciprocity laws'', such as Kato's explicit reciprocity law for the Beilinson--Kato elements \cite[Theorem 16.6]{kato04}, and the explicit reciprocity law of \cite[Theorem B]{KLZ17} for Beilinson--Flach elements. The conjectures of the preceding paragraphs suggest, at least to the present authors, that one should interpret any result comparing Euler systems of different ranks as an explicit reciprocity law.
  
 \section{Modular forms over an imaginary quadratic field}
  
  We now give an extended example showing some of the phenomena predicted by the conjectures of the previous sections. Many of the most interesting consequences only appear in situations where $\cK$ contains a $p$-adic Lie extension of dimension $> 1$; this can only occur, of course, when $K \ne \QQ$ (and, subject to Leopoldt's conjecture, if $K$ is not totally real). 
  
  We shall take $K$ to be an imaginary quadratic field with $p = \fp_1 \fp_2$ split in $K$, and suppose that $\cK$ includes the unique $\ZZ_p^2$-extension $F_\infty$ of $K$. We take $T = (T_f^*) |_{G_K}$, where $T_f$ is the representation of $G_{\QQ}$ associated to a weight 2 modular eigenform $f$. 
   
\begin{figure}[ht]
    \caption{Panchishkin subrepresentations for twists of $T$}
    \label{fig1}
    \begin{tikzpicture}
     [font=\footnotesize, 
      x=1.35cm, y=1.35cm, 
      baseline={([yshift={-2.5cm}]current bounding box.north)}
     ]
   
     \def\x{0.35}
     \colorlet{gy}{black!8!white} 
     
     \fill[fill=gy](-3.5, -1+\x) -- (-1, -1+\x) arc (90:0:\x) -- (-1+\x, -3.5) -- (-3.5, -3.5);
     \draw[thick](-3.5, -1+\x) -- (-1, -1+\x) arc (90:0:\x) -- (-1+\x, -3.5);
     
     \filldraw[thick, fill=gy](-3.5, -\x) -- (-1,  -\x) arc (-90:90:\x) -- (-1, \x) -- (-3.5, \x);
     
     \filldraw[thick, fill=gy](-\x, -3.5) -- (-\x, -1) arc (180:0:\x) -- (\x, -1) -- (\x, -3.5);
     
     \fill[fill=gy](-3.5, 1-\x) -- (-1, 1-\x) arc (-90:0:\x) -- (-1+\x, 3.5) -- (-3.5, 3.5);
     \draw[thick, fill=gy](-3.5, 1-\x) -- (-1, 1-\x) arc (-90:0:\x) -- (-1+\x, 3.5);
     \fill[fill=gy](-1+\x, 3.5) -- (-3.5, 3.5) -- (-3.5, 1-\x);
     
     \fill[fill=gy](3.5, -1+\x) -- (1, -1+\x) arc (90:180:\x) -- (1-\x, -3.5) -- (3.5, -3.5);
     \draw[thick](3.5,-1+\x) -- (1, -1+\x) arc (90:180:\x) -- (1-\x, -3.5);
     
     \filldraw[thick, fill=gy] (\x, 0) arc (0:360:\x);
     
     \foreach \x in {-3, ..., 3} {
      \foreach \y in{-3, ..., 3} {
       \fill[color=black] (\x, \y) circle (1pt);
      }
     }
     
     \draw[ semithick, ->] (-3.5, 0) -- (3.5, 0) node[anchor=west] {$a$};
     \draw[ semithick, ->] (0, -3.5) -- (0, 3.5) node[anchor=south] {$b$};
     
     \foreach \x in {1, ..., 3} {
      \draw[ semithick] (\x, 0.1) -- (\x, 0) node[anchor=north]{\footnotesize{\x}};
      \draw[ semithick] (-0.1, \x) -- (0, \x) node[anchor=west]{\footnotesize{\x}};
     }
     \foreach \x in {-1, ..., -2} {
      \draw[ semithick] (\x, 0.1) -- (\x, 0) node[anchor=north]{\footnotesize{\x}};
      \draw[ semithick] (-0.1, \x) -- (0, \x) node[anchor=west]{\footnotesize{\x}};
     }
     
     \node[draw=black, fill=white] () at (-3, -3) {$\Sigma^{(4)}$};
     \node[draw=black, fill=white] () at (-3, 3) {$\Sigma^{(2)}$};
     \node[draw=black, fill=white] () at (3, -3) {$\Sigma^{(2')}$};
     \node[draw=black, fill=white] () at (-3, 0) {\footnotesize{$\Sigma^{(3)}$}};
     \node[draw=black, fill=white] () at (0, -3) {\footnotesize{$\Sigma^{(3')}$}};
     \node[draw=black, fill=white, anchor=south west] () at (0.8*\x, 0.8*\x) {$\Sigma^{(1)}$};
   
    \end{tikzpicture}
    
    \bigskip

    \begin{tabular}{|c|c|c|c|}
     \hline
     Region & Critical? & $T_{\fp_1}^+$ & $T_{\fp_2}^+$ \\
     \hline
     \rule{0pt}{2.5ex} 
     $\Sigma^{(1)}$ & 0-crit & $T^+$ & $T^+$ \\
     $\Sigma^{(2)}$ & 0-crit & $T$   & 0 \\
     $\Sigma^{(2')}$ & 0-crit & 0   & $T$ \\
     \hline
     \rule{0pt}{2.5ex}
     $\Sigma^{(3)}$ \strut & 1-crit & $T^+$ & $T$ \\
     $\Sigma^{(3')}$ & 1-crit & $T$ & $T^+$ \\
     \hline
     \rule{0pt}{2.5ex}
     $\Sigma^{(4)}$ & 2-crit & $T$ & $T$ \\
     \hline
    \end{tabular}
   \end{figure}
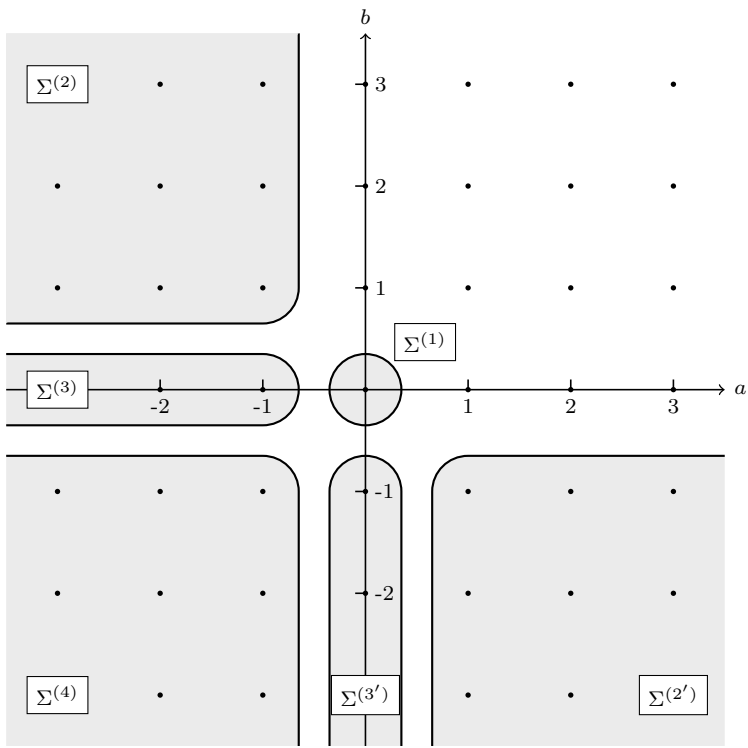

   If $\chi$ is a character of $\Gal(F_\infty / K)$ which is de Rham above $p$ (and hence corresponds to an algebraic Gr\"ossencharacter of $K$), then $\chi$ has two Hodge--Tate weights $(a,b)$. In Figure \vref{fig1} (adapted from Figure 1 of \cite{leiloefflerzerbes14b}), the shaded areas are the regions of the $(a, b)$ plane for which $T(\chi^{-1})$ is $r$-critical for some $r \ge 0$. Assuming $f$ is ordinary at $p$, so there is a 1-dimensional subrepresentation $T^+$ of $T_f^*$ at $p$, we can describe Panchishkin subrepresentations for each of these regions as in the accompanying table.
   
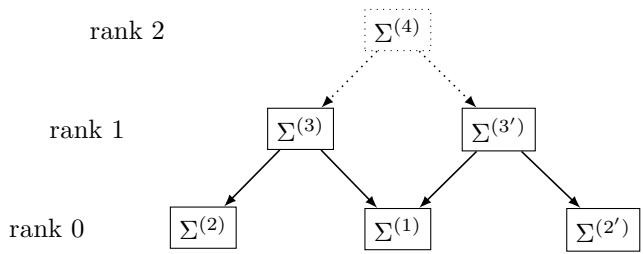
\begin{figure}[bt]
    \caption{Euler systems and explicit reciprocity laws for $T$}
    \label{fig2}
    \vspace{0.5cm}
    \begin{tikzpicture}[>=latex]
     \node[draw=black,dotted] (X) at (0,0) {$\Sigma^{(4)}$};
     \node[draw=black] (Y) [below left=0.75cm and 0.4cm of X]  {$\Sigma^{(3)}$};
     \node[draw=black] (Z) [below right=0.75cm and 0.4cm of X] {$\Sigma^{(3')}$};
     \node[draw=black] (U) [below left=0.75cm and 0.4cm of Y]  {$\Sigma^{(2)}$};
     \node[draw=black] (V) [below left=0.75cm and 0.4cm of Z]  {$\Sigma^{(1)}$};
     \node[draw=black] (W) [below right=0.75cm and 0.4cm of Z]  {$\Sigma^{(2')}$};
     
     \draw [semithick,dotted,->] (X) -- (Y);
     \draw [semithick,dotted,->] (X) -- (Z);
     \draw [semithick,->] (Y) -- (U);
     \draw [semithick,->] (Y) -- (V);
     \draw [semithick,->] (Z) -- (V);
     \draw [semithick,->] (Z) -- (W);
     
     \node () [left=2.5cm of X] {rank 2};
     \node () [left=1.75cm of Y] {rank 1};
     \node () [left=1cm of U] {rank 0};
    \end{tikzpicture}
   \end{figure}

   So Conjecture \ref{conj:prefinement} predicts that we should have six Euler systems in this setting, one for each region in the diagram: one of rank 2, two of rank 1, and three of rank 0. Moreover, these should be connected by explicit reciprocity laws corresponding to the the arrows in Figure \vref{fig2}.
  
   At present, the bottom half of Figure \ref{fig2} (the part drawn in solid ink) is firmly established. The three rank 0 Euler systems -- or at least their $p$-parts, which are measures on $\Gal(F_\infty / K)$ -- are familiar objects: they are the three $p$-adic $L$-functions described in \cite[Theorem 6.1.3]{leiloefflerzerbes14b}. The two rank 1 Euler systems can be constructed using Beilinson--Flach elements associated to CM families of modular forms; the construction of the CM family relies on a choice of prime above $p$, so one obtains two Euler systems corresponding to the regions $\Sigma^{(3)}$ and $\Sigma^{(3')}$. The four arrows linking these to the $p$-adic $L$-functions are all instances of the explicit reciprocity law of \cite[Theorem B]{KLZ17}. However, the top, dotted half of the diagram is more mysterious, since we know of no plausible geometric approach to constructing a rank 2 Euler system for the twists in $\Sigma^{(4)}$.
   
   \begin{remark} \
    \begin{enumerate}
     \item The $p$-adic $L$-function associated to $\Sigma^{(1)}$ can actually be defined over a finite extension of $\Qp$ (instead of the rather large, but still discretely-valued, extension $L$). However, those for $\Sigma^{(2)}$ and $\Sigma^{(2')}$ do not descend in any canonical way. More subtly, the base extension to $L$ is also needed in order to define the rank 1 Euler systems for $\Sigma^{(3)}$ and $\Sigma^{(3')}$: the Beilinson--Flach elements \emph{a priori} take values in $V_f^* \otimes V_\mathbf{g}^*$ where $\mathbf{g}$ is an auxiliary CM Hida family induced from $K$. To identify them with classes in $V_f^*$ alone, we need to find a basis of $V_\mathbf{g}^*$ in which $G_K$ acts diagonally. There is no canonical choice of such a basis over $\Qp$, but after base-extending to $L$ we can obtain a canonical basis from Ohta's $\Lambda$-adic comparison isomorphism.
     
     \item We can obtain a Panchishkin subrepresentation for twists in $\Sigma^{(1)}$ without assuming that $p$ is split, but the assumption that $f$ be ordinary is essential. On the other hand, for $\Sigma^{(2)}$ and its mirror-image $\Sigma^{(2')}$, the ordinarity of $f$ is not needed, but the splitting of $p$ is essential; and both conditions are needed simultaneously for $\Sigma^{(3)}$ or for $\Sigma^{(3')}$. These are, of course, special cases of the remarks about Rankin--Selberg convolutions in section \ref{sect:exRS} above, since the $L$-function $L(f / K, \chi, s)$ can also be described as the Ran\-kin--Selberg convolution of $f$ with a CM form induced from $\chi$.
     
     \item As sketched in \S \ref{sect:mainconj} above, for each node in Figure \ref{fig2} we can formulate an Iwasawa main conjecture of Greenberg type for $T$ over $F_\infty$. These conjectures are not independent of each other: an argument using Poitou--Tate duality shows that whenever two nodes are related by an explicit reciprocity law, the corresponding main conjectures are equivalent. It follows, for instance, that the Greenberg--Iwasawa main conjectures for $\Sigma^{(1)}$ and $\Sigma^{(2)}$ are equivalent to each other. Although there is no direct link between the $p$-adic $L$-functions concerned, they are tied together by the explicit reciprocity laws relating both of them to the rank 1 Euler system associated to $\Sigma^{(3)}$. This observation is due to Xin Wan, and its generalisations play a prominent role in recent work of Wan and his coauthors on the cyclotomic Iwasawa main conjecture and BSD leading term formula for supersingular elliptic curves over $\QQ$ \cite{wan15, jetchevskinnerwan17}.
     
     \item The representation $T_f |_{G_K}$ is the Galois representation attached to the base-change of $f$ to $K$, which is a cohomological automorphic form for the group $\GL_2/K$. The conjectural picture of Euler systems for $T$ that we describe here would apply equally to  the $G_K$-representation attached to any cohomological eigenform $\cF$ for $\GL_2 / K$, whether or not it arises from base-change, as long as $\cF$ is ordinary at $\fp_1$ and $\fp_2$. However, in the non-base-change setting we can prove much less; for instance, we know of no way of $p$-adically interpolating the values $L(\cF/K, \chi, 0)$ for $\chi \in \Sigma^{(2)}$ if $\cF$ is a non-base-change form.
    \end{enumerate}
    \end{remark}
   
  \section{The non-ordinary case}
  
  Greenberg's formulation of Iwasawa theory relies on the existence of Panchkishkin subrepresentations, but in many interesting cases these do not exist. A more flexible theory has been developed by Pottharst \cite{pottharst13}, based on the observation that for each $v \mid p$ one can attach to $V |_{G_{K_v}}$ a semilinear algebra object known as a $(\varphi, \Gamma)$-module, denoted $D^{\dag}_{\mathrm{rig}}(K_v, V)$; and there may be interesting subobjects of $D^{\dag}_{\mathrm{rig}}(K_v, V)$ which do not come from subrepresentations of $V$. For instance, if $f$ is a modular form, one may attach a rank-1 submodule of $D^{\dag}_{\mathrm{rig}}(\Qp,V_f)$ to any non-zero root $\alpha$ of the Hecke polynomial $X^2 - a_p(f) X + p^{k-1} \varepsilon_f(p)$, while this submodule only comes from a subrepresentation if $\alpha$ is a $p$-adic unit.
  
  The downside of working with these objects is that one has to give away some control of denominators: the ``analytic Iwasawa cohomology'' modules appearing in Pottharst's theory are not modules over the Iwasawa algebra $\Lambda(\Gamma)$, but over the larger algebra $\mathcal{H}(\Gamma)$ of locally analytic distributions on $\Gamma$, which is a $\Qp$-algebra having no natural $\Zp$-lattice. So, in translating from the classical language to the new one, we lose control of the $\mu$-invariants of Selmer groups. 
  
  Subject to this caveat, one can generalise the entire conjectural picture of Euler systems described above assuming only that one has a ``Panchishkin submodule'' of $D^{\dag}_{\mathrm{rig}}(K_v, V)$ for each $v \mid p$, i.e.~a subobject which precisely accounts for all the positive Hodge--Tate weights. When this occurs, we should expect to be able to extend the cohomology classes of Conjecture \ref{conj:basicconjecture} to elements of Pottharst's analytic cohomology modules in the $p$-direction, and these Euler systems should satisfy main conjectures, formulated in terms of equalities of characteristic ideals over $\mathcal{H}(\Gamma)$.
  
  \begin{remark}
   One new phenomenon that occurs when one recasts the theory in Pottharst's setting is that Panchishkin submodules are no longer unique. Hence one should formulate Conjecture \ref{conj:pdirection} as associating a family of elements of the $r$-th powers of Pottharst's analytic cohomology modules to an $r$-critical $G_K$-representation together with a \emph{choice} of Panchishkin submodule at each $v \mid p$ (which should be understood as a ``$p$-stabilisation''). For instance, the non-ordinary analogue of Figure \ref{fig2} consists of 11 objects (one Euler system of rank 2, four of rank 1, and six of rank 0), with 16 potential explicit reciprocity laws connecting them. Many, but not all, of these can be constructed using the techniques of \cite{loefflerzerbes16}.
  \end{remark}

  \providecommand{\bysame}{\leavevmode\hbox to3em{\hrulefill}\thinspace}
  \providecommand{\MR}[1]{\relax}
  \renewcommand{\MR}[1]{%
   MR \href{http://www.ams.org/mathscinet-getitem?mr=#1}{#1}.
  }
  \providecommand{\href}[2]{#2}
  \newcommand{\articlehref}[2]{\href{#1}{#2}}

\end{document}